\documentclass[a4paper,12pt]{article}

\usepackage{amsmath,amssymb,amsthm}

\newcommand{\A}{{\mathcal A}}
\newcommand{\Z}{{\mathbb Z}}
\newcommand{\R}{{\mathbb R}}
\newcommand{\C}{{\mathbb C}}

\newcommand{\h}{{\mathfrak h}}
\newcommand{\Lyndon}{{\mathfrak L}}

\newcommand{\sing}{\Sigma}
\newcommand{\ncpa}[1]{\C\langle#1\rangle}
\newcommand{\ncfps}[1]{\ncpa{\langle#1\rangle}}
\newcommand{\reg}{\mathrm{reg}}
\newcommand{\Liz}[4]{\mathrm{Li}_{#1}^{#2}\left(#3;#4\right)}
\newcommand{\Li}[2]{\mathrm{Li}_{#1}^{#2}}
\newcommand{\z}[2]{\left(#1;#2\right)}
\newcommand{\imu}{\sqrt{-1}}
\newcommand{\Lv}{{\mathcal L}}

\DeclareFontEncoding{OT2}{}{}
\DeclareFontSubstitution{OT2}{cmr}{m}{n}
\DeclareSymbolFont{cyss}{OT2}{wncyss}{m}{n}
\DeclareMathSymbol{\sh}{\mathbin}{cyss}{`x}

\newtheorem{thm}{Theorem}
\newtheorem{prop}[thm]{Proposition}
\theoremstyle{definition}
\newtheorem{exam}{Example}

\title
{Duality formulas of the Special Values of Multiple Polylogarithms}
\author{OKUDA Jun-ichi}
\date{}

\usepackage[dvips]{hyperref}

\begin{document}

\maketitle

\begin{abstract}
 The special values of multiple polylogarithms,
 which include multiple zeta values,
 appear in several fields of mathematics and physics.
 Many kinds of their linear relations are investigated
 as well as their algebraic relations.
 From the viewpoint of a connection matrix of Fuchsian equations,
 two kinds of duality of these values are derived.
\end{abstract}

\section{Introduction}

In \cite{MR1709772} the multiple polylogarithm is defined by
\begin{align*}
  \lambda
  \begin{pmatrix}
   s_1&\dots& s_k\\
   b_1&\dots& b_k
  \end{pmatrix}
 =
 \sum_{m_1=k}^\infty
 \sum_{m_1>m_2>\dots>m_k>0}
 \frac{b_1^{-(m_1-m_2)} \dotsb b_{k-1}^{-(m_{k-1}-m_k)} b_k^{-m_k}}
  {m_1^{s_1}\dotsb m_{k-1}^{s_{k-1}} m_k^{s_k}},
\end{align*}
and many relations, reductions, and explicit evaluations
of special values of them are investigated.
In \cite{AK2002}
a class of the special values of the multiple polylogarithm
are called multiple $L$-values.
For any $m\in\Z_{>0}$,
we put $R_m=\Z \big/ m\Z$ and $\zeta_m = e^{2\pi\imu/m}$.
The multiple $L$-value of positive integers $k_i$ 
and $a_i\in R_m$ ($i=1,\dots,r$)
is defined by
\begin{align*}
 L(k_1,\dots,k_r;a_1,\dots,a_r)
 =
 \lambda
 \begin{pmatrix}
  k_1& \dots& k_{r-1}& k_r\\
  \zeta_m^{-a_1} & \dots& \zeta_m^{-a_{r-1}}& \zeta_m^{-a_r}
 \end{pmatrix}.
\end{align*}
In particular the case $m=1$ is known as multiple zeta values
\begin{align*}
 \zeta(k_1,k_2,\dots,k_r)
 =
 \sum_{m_1>m_2>\dots>m_r>0}
 \frac{1}{m_1^{k_2}m_2^{k_2}\dotsm m_r^{k_r}},
\end{align*}
and many kinds of relations and properties are discovered.
Also the case $m=6$ is discussed in \cite{MR2002a:81180}
from the point of view of quantum field theory
and some relations are obtained.
The main purpose of the study of these values is
to discover all the linear relations explicitly.

In \cite{MR92f:16047} the KZ associator $\Phi(X,Y)$ is defined by
\begin{align*}
 \Phi(X,Y) = G_1^{-1}(z) G_0(z),
\end{align*}
where $G_0$ and $G_1$ are the solutions of the Fuchsian equation
called the formal KZ equation
\begin{align}
 \frac{dG}{dz}
  = 
    \left(\frac{X}{z} + \frac{Y}{z-1}\right) G
 ,
 \qquad G(z) \in \ncfps{X,Y},
 \label{eq:formal KZ}
\end{align}
characterized by the following asymptotic properties respectively:
\begin{align*}
 G_0(z) \times z^{-X} &\to 1 \qquad (z\to 0)&
 &\text{and}&
 G_1(z) \times (1-z)^{-Y} &\to 1 \qquad (z\to 1).
\end{align*}
Here $z^{-X}=\exp\left(\sum_{n=0}^\infty (-X\log z)^n/n!\right)$
and so on.
It is conjectured that
all the relations for multiple zeta values
can be deduced from the relations which the KZ associator satisfies.
So the KZ associator is a very important generating function
of multiple zeta values.
For example, the coefficients of the relation
\begin{align}
 \Phi(X,Y)^{-1} = \Phi(Y,X)
 \label{eq: associator dual}
\end{align}
give rise to the duality formula for multiple zeta values.
From this point of view some studies have been developed
(c.f.\ \cite{MPH1998,Fur2003,math.NT/0301277}).

In this paper we consider the next generalization of the KZ associator.
Using the connection matrices of specified solutions,
we investigate the relations among
the special values of multiple polylogarithms
including multiple $L$-values.
For any finite subset  $\sing\subset \C$
we define the generalization of the formal KZ equation as follows:
\begin{align}
 \frac{dH}{dz}
 =
 \left(\sum_{a\in\sing} \frac{X_{a}}{z-a}\right) H,
 \qquad
 H(z)\in\ncfps{X_{a};a\in\sing}.
 \label{eq:ext KZ}
\end{align}
This differential equation is also Fuchsian with singular points
$\sing\cup\{\infty\}$,
so for any $a,b\in\sing$
there is a unique local solution $H^{ab}_\sing(z)$ 
with the asymptotic property
\begin{align*}
 H^{ab}_\sing(z) \times \left(\frac{z-b}{a-b}\right)^{-X_b}
 &\longrightarrow 1
 &(z&\to b).
\end{align*}
Then we define the connection matrix
$\Phi^{ab}_\sing(X_c;c\in\sing)\in\ncfps{X_c;c\in\sing}$ by
\begin{align*}
 \Phi^{ab}_\sing(X_c;c\in\sing) = H^{ba}_\sing(z)^{-1}  H^{ab}_\sing(z).
\end{align*}
This connection matrix is the counterpart of the KZ associator.

This article is organized as follows:
In Section~\ref{sec: shuffle algebra}
the shuffle algebra w.r.t.\ \eqref{eq:ext KZ} is defined
and some properties are described.
In Section~\ref{sec: multilog}
by using the shuffle algebra,
we define the multiple polylogarithms
and construct the solutions of \eqref{eq:ext KZ}.
Comparing the two evaluations of the ratio of solutions,
we obtain the linear relations of special values of multiple polylogarithms.
In Section~\ref{sec: dual}
we derive the dual and closed formula of the linear relations
of special values of multiple polylogarithms
from the symmetry of singular points of \eqref{eq:ext KZ}.

\section{Shuffle Algebra}\label{sec: shuffle algebra}
Following \cite{IKZ2003,AK2002},
we introduce the shuffle algebra and its properties.
We define a non-commutative polynomial algebra
$\A_{\sing} := \ncpa{x_{a}; a\in \sing}$,
and call the elements $\{x_a;a\in\sing\}$ letters
     and the monomials in $\A_{\sing}$ words.
The weight of the word is defined as the total number of letters
which appear in it.
In this algebra we define the shuffle product ``$\sh$''
recursively as follows:
\begin{enumerate}
 \item $w\sh 1 = 1\sh w = w$,
 \item $l_1w_1\sh l_2w_2 = l_1(w_1\sh l_2w_2) + l_2(l_1w_1\sh w_2)$,
\end{enumerate}
where $l_1$ and $l_2$ are letters, and $w_1, w_2$ and $w$ are words.
We regard $\A_\sing = (\A_\sing, +, \sh)$ and call it
the ``shuffle algebra''.
For example the $n$-th power of a letter $x_a$ is computed as
\begin{align}
 x_a^{\sh n}
 = \underbrace{x_a\sh \dots\sh x_a}_{n}
 = n! x_a^n.
 \label{eq:n-th power}
\end{align}

Fix any two elements $a$ and $b\in\sing$
and let ``$\prec$'' be a total order of $\sing$
with the maximum element $a$ and the minimum element $b$.
From this order, the words in $\A_\sing$ are totally ordered
by the lexicographic order.
The Lyndon words $\Lyndon=\Lyndon_{\sing,\prec}\subset\A_\sing$ are defined as follows:
\begin{align*}
 \Lyndon = \{w\neq1\in\A_\sing\text{: word} \mid
 \text{for any non-trivial decomposition } w=uv,\ w \prec v\}.
\end{align*}
In particular $\Lyndon\supset\{x_a\}_{a\in\sing}$
and
Lyndon words which start from $x_a$ (respectively ends at $x_b$)
are only $x_a$ (respectively $x_b$).
Then the shuffle algebra $\A_\sing$ is the commutative polynomial algebra
$\C[\Lyndon]$ (\cite{MR94j:17002} Theorem 6.1.).
We define subalgebras of $\A_\sing$ as
\begin{align*}
 \A^b_\sing
  &:=
  \C[\Lyndon\backslash\{x_b\}]
  =
  \C.1 \oplus \bigoplus_{\substack{l\in\sing\\l\neq b}}\A_\sing x_l,\\
 \A^{ab}_\sing
  &:=
  \C[\Lyndon\backslash\{x_a,x_b\}]
  =
  \C.1
  \oplus \bigoplus_{\substack{l_1,l_2\in\sing\\l_1\neq a,l_2\neq b}}
     x_{l_1} \A_\sing x_{l_2}.
\end{align*}
Because $x_a$ and $x_b \in L$,
$\A_\sing$ can be written as follows:
\begin{align}
 \A_\sing &= \A^b_\sing[x_b]
          := \bigoplus_{j=0}^\infty \A_\sing^b\sh x_b^{\sh j}
            \label{eq:b poly}\\
          &= \A^{ab}_\sing[x_a, x_b]
          := \bigoplus_{i,j=0}^\infty x_a^{\sh i}\sh\A_\sing^{ab}\sh x_b^{\sh j}.
\end{align}
Using these decompositions,
we define $\reg^b: \A_\sing \to \A_\sing^b$
and $\reg^{ab}: \A_\sing \to \A_\sing^{ab}$ to be the maps
sending any word to its constant term of $\A_\sing^{b}[x_b]$ and
$\A_\sing^{ab}[x_a,x_b]$ respectively.
By definition, $\reg^b$ and $\reg^{ab}$ are $\sh$-homomorphisms.

\begin{prop}[\cite{IKZ2003}]
 For $w=w_bx_b^n=x_a^mw_{ab}x_b^n\in\A_\sing$
 ($w_b\in\A_\sing^{b}, w_{ab}\in\A^{ab}_\sing$)
 we have
 \begin{align}
  &\reg^b(w_b x_b^n) = \sum_{j=0}^n (-1)^j w_bx_b^{n-j} \sh x_b^j,
   \label{eq:b reg}\\
  &\reg^{ab}(x_a^m w_{ab} x_b^n)
   = \sum_{i=0}^m\sum_{j=0}^n (-1)^{i+j} x_a^i \sh x_a^{m-i}w_bx_b^{n-j} \sh x_b^j,\\
  &w = \sum_{j=0}^n \reg^b(w_b x_b^{n-j})\sh x_b^j
    = \sum_{i=0}^m \sum_{j=0}^n x_a^i \sh \reg^{ab}(x_a^{m-i}w_{ab} x_b^{n-j})\sh x_b^j.
  \label{eq:w decomp}
 \end{align}
\end{prop}

Prepare the generating function of all words of $\A_\sing$ in the non-commutative
power series algebra $\A_\sing\langle\langle X_a;a\in\sing\rangle\rangle$
\begin{align*}
 \sum_W wW
  &=
 \begin{aligned}[t]
    &1 + x_{a_0}X_{a_0} + x_{a_1}X_{a_1} + x_{a_2}X_{a_2}\dotsb\\
   &+ x_{a_0}x_{a_0}X_{a_0}X_{a_0} + x_{a_0}x_{a_1}X_{a_0}X_{a_1}
    + x_{a_0}x_{a_2}X_{a_0}X_{a_2} + \dotsb\\
   &\ + x_{a_1}x_{a_0}X_{a_1}X_{a_0} + x_{a_1}x_{a_1}X_{a_1}X_{a_1}
    + x_{a_1}x_{a_2}X_{a_1}X_{a_2} + \dotsb,
 \end{aligned}
\end{align*}
where the sum is taken over all words
in $\A_\sing\langle\langle X_a;a\in\sing\rangle\rangle$
and $W$ is the capitalization of $w$.
Noting \eqref{eq:n-th power} and using this generating function,
\eqref{eq:w decomp} can be expressed by
\begin{align}
\sum_W wW
 &=
 \left( \sum_W \reg^{b}(w) W\right)
 \times
 \exp\left(x_bX_b \right) \label{eq:b decomp}\\
 &=
 \exp\left(x_aX_a \right)
 \times
 \left( \sum_W \reg^{ab}(w) W\right)
 \times
 \exp\left(x_bX_b \right), \label{eq:ab decomp}
\end{align}
where $\exp\left(xX \right)$ means
\begin{align*}
 \exp\left(xX \right)
 =
 \sum_{n=0}^\infty
 \frac{x^{\sh n}}{n!} X^n
 =
 \sum_{n=0}^\infty
  x^n X^n.
\end{align*}

The inverse $\left(\sum_W wW\right)^{-1}$ is
\begin{align*}
 \left(\sum_W wW\right)^{-1} = \sum_W S(w)W = \sum_W wS(W)
\end{align*}
where $S$ is an anti-involution w.r.t.\ the ordinary product defined by
$S: x_a \mapsto -x_a, X_a \mapsto -X_a$.
In general any homomorphism or anti-homomorphism
w.r.t.\ ordinary product is a homomorphism w.r.t.\ $\sh$.
So $S$ is also a $\sh$-homomorphism.

\section{Multiple Polylogarithms}\label{sec: multilog}
For a word
$w= x_b^{k_1-1}x_{c_1}\dotsm x_b^{k_{r-1}-1}x_{c_{r-1}}x_b^{k_r-1}x_{c_r}
\in\A_\sing^b$
($b\neq c_i\in\sing$)
we define the multiple polylogarithm (of one variable)
$\Liz{\sing}{b}{w}{z}$ by
\begin{align*}
 \Liz{\sing}{b}{1}{z} &:= 1,\\
 \Liz{\sing}{b}{w}{z}
 &:=
 (-1)^r
 \lambda
 \begin{pmatrix}
  k_1& \dots& k_{r-1}& k_r\\
  \frac{c_1-b}{z-b}& \dots& \frac{c_{r-1}-b}{z-b}& \frac{c_r-b}{z-b}
 \end{pmatrix}\\
 &=
 (-1)^r
 \sum_{m_1>\dots>m_r>0}
 \frac{
  \left(\frac{z-b}{c_1-b}\right)^{m_1-m_2}
  \dotsb
 \left(\frac{z-b}{c_{r-1}-b}\right)^{m_{r-1}-m_{r}}
  \left(\frac{z-b}{c_r-b}\right)^{m_{r}}
 }
 {m_1^{k_1}\dotsb m_{r-1}^{k_{r-1}}m_r^{k_r}},
\end{align*}
where the series converges absolutely
for $|z-b|<\min_{c\in\sing}|c-b|$.
If the weight of $w$ is greater than or equal to $1$,
$\Liz{\sing}{b}{w}{b}=0$.
Moreover Li is extended linearly w.r.t.\ $w$.
For example
\begin{align}
 \Liz{\sing}{b}{x_c}{z}
 = (-1)\sum_{m=1}^\infty
   \frac{ \left( \frac{z-b}{c-b} \right)^m }{m}
 = \log\left(1-\frac{z-b}{c-b}\right)
 = \log\left(\frac{z-c}{b-c}\right).
 \label{eq: exam Li}
\end{align}
$\Liz{\sing}{b}{w}{z}$ is easily expressed as the iterated integral
\begin{multline*}
  \Liz{\sing}{b}{w}{z}=
  \underbrace{\int_{b}^z \frac{dz}{z-b} \dotsb \int_b^{z}
  \frac{dz}{z-b} \int_b^z \frac{dz}{z-c_1} }_{k_1}\dotsm\\
  \dotsm
  \underbrace{\int_{b}^z \frac{dz}{z-b} \dotsb \int_b^{z}
  \frac{dz}{z-b} \int_b^z \frac{dz}{z-c_{r-1}} }_{k_r-1}
  \underbrace{\int_{b}^z \frac{dz}{z-b} \dotsb \int_b^{z}
  \frac{dz}{z-b} \int_b^z \frac{dz}{z-c_r} }_{k_r}.
\end{multline*}
From this expression $\Liz{\sing}{b}{w}{z}$ can be analytically continued
to $\C\backslash\sing$ and the derivative is
\begin{align}
 \frac{d}{dz}\Liz{\sing}{b}{x_cw}{z} = \frac{1}{z-c}\Liz{\sing}{b}{w}{z}
 \label{eq:diff rel}
\end{align}
for any letter $c$ and any word $w\in\A_\sing^b$.
By using an induction w.r.t. the sum of weight of words,
we obtain the next proposition.

\begin{prop}
 For any words $w_1, w_2 \in \A_\sing^b$,
 \begin{align*}
  \Liz{\sing}{b}{w_1}{z}\Liz{\sing}{b}{w_2}{z} = \Liz{\sing}{b}{w_1\sh w_2}{z}.
 \end{align*}
\end{prop}

This proposition inspires us to define $\Liz{\sing}{ab}{w}{z}$
for any word $w\in\A_\sing$ by
\begin{align}
 \Liz{\sing}{ab}{w}{z} & := \Liz{\sing}{b}{w}{z}
  \quad \text{for $w\in\A_\sing^b$},\\
 \Liz{\sing}{ab}{x_b}{z} &:= \log\left(\frac{z-b}{a-b}\right)
 = \Liz{\sing}{a}{x_b}{z},
 \label{eq: def Li}
\end{align}
and extend as a $\sh$-homomorphism w.r.t.\ $\A_\sing$.
It is well-defined because of \eqref{eq:b poly}.

\begin{thm}
 \begin{align*}
  H^{ab}_\sing(z) = \sum_W \Liz{\sing}{ab}{w}{z} W.
 \end{align*}
\end{thm}
\begin{proof}
 Applying $\Liz{\sing}{ab}{\cdot}{z}$ to \eqref{eq:b decomp} we have
 \begin{align*}
  \sum_W \Liz\sing {ab}wz W
  &= \left(\sum_W \Liz{\sing}{b}{\reg^b(w)}{z}W\right)
    \times
    \exp \left(\Liz{\sing}{ab}{x_b}{z} X_b\right)\\
  &=
   \left(\sum_W \Liz{\sing}{b}{\reg^b(w)}{z}W\right)
    \times
    \left(\frac{z-b}{a-b}\right)^{X_b}.
 \end{align*}
 So the asymptotic property is satisfied.
 We must show that this series is the solution of \eqref{eq:ext KZ},
 i.e.\ \eqref{eq:diff rel} holds for any $w\in\A_\sing$.
 Applying \eqref{eq:b reg} to the word $x_c w_b x_b^n$
 ($w_b\in\A_\sing^b$) we have
 \begin{align}
  \reg^b(x_c w_b x_b^n)
  &=
  \sum_{j=0}^n (-1)^j x_c w_b x_b^{n-j} \sh x_b^{j}\\
  &=
  x_c\sum_{j=0}^n (-1)^j w_b x_b^{n-j} \sh x_b^{j}
  + x_b \sum_{j=1}^n (-1)^j x_c w_b x_b^{n-j} \sh x_b^{j-1}\\
  &= x_c \reg^b(w_b x_b^n) - x_b \reg^b(x_c w_b x_b^{n-1}).
  \label{eq: reg part}
 \end{align}
 With this formula and \eqref{eq:w decomp}, we can calculate
 the derivative of $\Liz{\sing}{ab}{x_c w_b x_b^n}{z}$ as follows:
 \begin{align*}
  &\frac{d}{dz}\Liz{\sing}{ab}{x_cw_b x_b^n}{z}\\
  &=
  \frac{d}{dz}
   \sum_{j=0}^n \Liz{\sing}{ab}{\reg^b\left(x_cw_b x_b^{n-j}\right)}{z}
                     \Liz{\sing}{ab}{x_b^j}{z}\\
  \begin{split}
   &=
  \frac{1}{z-c}\Liz{\sing}{ab}{\reg^b\left(w_b x_b^{n}\right)}{z}
  -\frac{1}{z-b}\Liz{\sing}{ab}{\reg^b\left(w_b x_b^{n-1}\right)}{z}\\
  &\quad
  +\sum_{j=1}^{n-1}
  \begin{aligned}[t]
  &\left\{
    \frac{1}{z-c}\Liz{\sing}{ab}{\reg^b\left(w_b x_b^{n-j}\right)}{z}
     \Liz{\sing}{ab}{x_b^j}{z}
   \right.\\
  &\quad
    -\frac{1}{z-b}\Liz{\sing}{ab}{\reg^b\left(x_cw_b x_b^{n-j-1}\right)}{z}
     \Liz{\sing}{ab}{x_b^j}{z}
   \\
  &\qquad+\left.
   \frac{1}{z-b}\Liz{\sing}{ab}{\reg^b\left(x_cw_b x_b^{n-j}\right)}{z}
    \Liz{\sing}{ab}{x_b^{j-1}}{z}
  \right\}
  \end{aligned}\\
  &\quad+
   \begin{aligned}[t]
  &\frac{1}{z-c} \Liz{\sing}{ab}{\reg^b(w_b)}{z} \Liz{\sing}{ab}{x_b^n}{z}\\
  &\qquad
  +\frac{1}{z-b} \Liz{\sing}{ab}{\reg^b(x_cw_b)}{z} \Liz{\sing}{ab}{x_b^{n-1}}{z}
   \end{aligned}  
  \end{split}\\
  &=
  \frac{1}{z-c}
  \sum_{j=0}^{n}
   \Liz{\sing}{ab}{\reg^b\left(w_b x_b^{n-j}\right)}{z}
    \Liz{\sing}{ab}{x_b^j}{z}
  \\
  &=
   \frac{1}{z-c} \Liz{\sing}{ab}{w_b x_b^n}{z}.
\end{align*}
\end{proof}

Using  \eqref{eq:ab decomp}, \eqref{eq: exam Li} and \eqref{eq: def Li}
$H^{ab}_\sing(z)$ also can be written as
\begin{align*}
 H^{ab}_\sing(z)
 =
 \left( \frac{z-a}{b-a} \right)^{X_a}
 \left( \sum_W \Liz{\sing}{ab}{\reg^{ab}(w)}{z} W\right)
 \left(\frac{z-b}{a-b}\right)^{X_b},
\end{align*}
and in the same way $H^{ba}_\sing(z)$ is
\begin{align*}
 H^{ba}_\sing(z)
 =
 \left( \frac{z-b}{a-b} \right)^{X_b}
 \left( \sum_W \Liz{\sing}{ba}{\reg^{ba}(w)}{z} W\right)
 \left(\frac{z-a}{b-a}\right)^{X_a}.
\end{align*}

Let $a$ be one of the nearest points of $\sing$ to $b$.
Then the ratio of these two solutions is
\begin{align}
\begin{split}
 &H^{ba}_\sing(z)^{-1}H^{ab}_\sing(z)\\
 &\quad=
 \begin{aligned}[t]
  \left(\sum_W\Liz{\sing}{ba}{S(w)}{z}W \right)
  \left(\sum_W\Liz{\sing}{ab}{w}{z}W \right)
 \end{aligned}
\end{split}\\
&\quad=
 \sum_W \left(
  \sum_{w_1w_2=w}
   \Liz{\sing}{ba}{S(w_1)}{z} \Liz{\sing}{ab}{w_2}{z}
 \right) W
 \label{eq:euler}
 \\
 &\quad=
 \begin{aligned}[t]
 &\left( \frac{z-a}{b-a} \right)^{-X_a}
  \left( \sum_W \Liz{\sing}{ba}{\reg^{ba}\circ S(w)}{z} W\right)
  \left( \frac{z-b}{a-b} \right)^{-X_b}\\
 &\quad\times
  \left( \frac{z-a}{b-a} \right)^{X_a}
  \left( \sum_W \Liz{\sing}{ab}{\reg^{ab}(w)}{z} W\right)
  \left( \frac{z-b}{a-b} \right)^{X_b}.
 \end{aligned}
 \label{eq:ratio}
\end{align}
On the other hand,
for any $w = x_b^{k_1-1}x_{c_1}\dotsm x_b^{k_r-1}x_{c_r}\in \A^{ab}$,
$\Liz\sing {ab} w z$ can be evaluated at $z=a$.
We define $\Lv^{ab}_\sing:\A_\sing^{ab}\to\C$ to be this evaluation:
\begin{align*}
 &\Lv^{ab}_\sing(w) := \lim_{t\to 1}\Liz\sing {ab} w {ta + (1-t)b}\\
 &=
 (-1)^r\sum_{m_1=r}^\infty \sum_{m_1>m_2>\dots>m_r>0}
  \frac{ \left(\frac{a-b}{c_1-b}\right)^{m_1-m_2}
          \dotsb
           \left(\frac{a-b}{c_{r-1}-b}\right)^{m_{r-1}-m_r}
            \left(\frac{a-b}{c_{r}-b}\right)^{m_r}
       }{m_1^{k_1}\dotsb m_{r-1}^{k_{r-1}} m_r^{k_r}}\\
 &=
  (-1)^r\lambda
  \begin{pmatrix}
   k_1&\dots& k_r\\
   \frac{c_1-b}{a-b}&\dots& \frac{c_r-b}{a-b}
  \end{pmatrix}
 .
\end{align*}
$\Lv_{\sing}^{ab}$ is obviously a $\sh$-homomorphism.
Tending $z$ to $a$ in \eqref{eq:ratio} we have
\begin{align*}
 \Phi_\sing^{ab}(X_c;c\in\sing)
  =  \sum_W \Lv_{\sing}^{ab}\left({\reg^{ab}(w)}\right) W.
\end{align*}
The coefficients of \eqref{eq:euler} yield the functional equation
\begin{align*}
 \sum_{w_1w_2=w}
   \Liz{\sing}{ba}{S(w_1)}{z} \Liz{\sing}{ab}{w_2}{z}
 =
 \Lv_{\sing}^{ab}\left({\reg^{ab}(w)}\right),
\end{align*}
for any word $w\in\A_\sing$.
We call this functional equation ``Euler's inversion formula''.
If $a$ and $b$ are the nearest points to each other,
we can also evaluate the ratio \eqref{eq:ratio} at $b$
and we have
\begin{align*}
 \Phi_\sing^{ab}(X_c;c\in\sing)
  =  \sum_W \Lv_{\sing}^{ba}\left({\reg^{ba}\circ S(w)}\right) W.
\end{align*}
Consequently for any word $w\in\A_\sing$ we obtain
\begin{align*}
 \Lv_{\sing}^{ab}\left({\reg^{ab}(w)}\right)
 =
 \Lv_{\sing}^{ab}\left({\reg^{ba}(S(w))}\right).
\end{align*}
We call this equation the ``duality formula''.

\begin{thm}
 For any $b\in\sing$ and
 one of the nearest points $a\in\sing$ to $b$,
 \begin{align*}
  \Phi_\sing^{ab}(X_c;c\in\sing)
  =
  \sum_W \Lv_{\sing}^{ab}\left({\reg^{ab}(w)}\right) W.
 \end{align*}
 In particular for any word $w\in\A_\sing^{ab}$ we have
 \begin{align*}
  \sum_{w_1w_2=w} \Liz{\sing}{ba}{S(w_1)}{z} \Liz{\sing}{ab}{w_2}{z}
  =
  \Lv_\sing^{ab}(w).
 \end{align*}
 Moreover, if $b$ is one of the nearest points to $a$, we have
 \begin{align}
  \Lv_\sing^{ba}(S(w)) = \Lv_\sing^{ab}(w).
  \label{eq:dual}
 \end{align}
\end{thm}

\begin{exam}
 Set the singular points to be $\sing=\{0,1\},a=1$ and $b=0$.
 Then for positive integers $a_i,b_i$ ($i=1,\dotsc,s$)
 the coefficients of the connection matrix are written as follows:
 \begin{align*}
 \begin{split}
  &\Lv_{\{0,1\}}^{10}
  ({x_0^{a_1}x_1^{b_1}x_0^{a_2}x_1^{b_2}\dotsb x_0^{a_s}x_1^{b_s}})\\
 &\quad=
 (-)^{\sum_{i=1}^s b_i}
 \zeta(\underbrace{a_1+1,1,\dotsc,1}_{b_1},
       \underbrace{a_2+1,1,\dotsc,1}_{b_2},
       \dotsc,
       \underbrace{a_s+1,1,\dotsc,1}_{b_s}),
 \end{split}\\
 \begin{split}
  &\Lv_{\{0,1\}}^{01}
  \left(
    S({x_0^{a_1}x_1^{b_1}x_0^{a_2}x_1^{b_2}\dotsb x_0^{a_s}x_1^{b_s}})
  \right)\\
  &\quad=
  (-)^{\sum_{i=1}^s a_i+b_i}
  \Lv_{\{0,1\}}^{01}
  \left(
    {x_1^{b_s}x_0^{a_s}\dotsb x_1^{b_2}x_0^{a_2}x_1^{b_1}x_0^{a_1}}
  \right)
 \end{split}\\
 &\quad=
 (-)^{\sum_{i=1}^s b_i}
 \zeta(\underbrace{b_s+1,1,\dotsc,1}_{a_s},
       \dotsc,
       \underbrace{b_2+1,1,\dotsc,1}_{a_2},
       \underbrace{b_1+1,1,\dotsc,1}_{a_1}).
 \end{align*}
 From \eqref{eq:dual} we have the equation
 \begin{multline*}
  \zeta(\underbrace{a_1+1,1,\dotsc,1}_{b_1},
       \underbrace{a_2+1,1,\dotsc,1}_{b_2},
       \dotsc,
       \underbrace{a_s+1,1,\dotsc,1}_{b_s})\\
 =
 \zeta(\underbrace{b_s+1,1,\dotsc,1}_{a_s},
       \dotsc,
       \underbrace{b_2+1,1,\dotsc,1}_{a_2},
       \underbrace{b_1+1,1,\dotsc,1}_{a_1}),
 \end{multline*}
 this is the duality formula for multiple zeta values (\cite{MR96k:11110}).

 Euler's inversion formula for $x_0x_1$ is
 \begin{multline*}
  \Liz{\{0,1\}}{01}{x_1x_0}{z}
  + \Liz{\{0,1\}}{01}{x_0}{z}  \Liz{\{0,1\}}{10}{x_1}{z}\\
  + \Liz{\{0,1\}}{10}{x_0x_1}{z}
  = \Lv_{\{0,1\}}^{10}(x_0x_1),
 \end{multline*}
 or equivalently
 \begin{align*}
  \mathrm{Li}_2(1-z)
  + \log(z)\log(1-z)
  + \mathrm{Li}_2(z)
  = \zeta(2).
 \end{align*}
 Here $\mathrm{Li}_2(z)=\sum_{n=1}^\infty z^n/n^2$
 is Euler's dilogarithm.
 This formula is well-known as Euler's inversion formula
 for the dilogarithm.
\end{exam}

\begin{exam}
 Set $\sing=\{0,1,-1\}, a=1$ and $b=0$.
 Any word $w_{10}\in\A_{\{0,1,-1\}}^{10}$ can be written as
 \begin{align*}
  w_{10}&=x_0^{k_1-1}x_{d_1}x_0^{k_2-1}x_{d_2}\dotsb
   x_0^{k_{r-1}-1}x_{d_{r-1}}x_0^{k_r-1}x_{d_r}
  \qquad
  (d_i = \pm 1),\\
\intertext{and $S(w_{01})$ also can be written as}
  S(w_{10})
  &=
  (-)^{\sum_{i=1}^{r'} k'_i}
  x_1^{k_1'-1}x_{e_1} x_2^{k_2'-1}x_{e_2} \dotsb
  x_1^{k_{r'-1}'-1} x_{e_{r'-1}}
  x_1^{k_{r'}'-1} x_{e_{r'}}
  \ (e_j=0,-1).
 \end{align*} 
 Then the duality formula gives us the equation
 \begin{multline}
   (-)^{|w_{01}|_{x_{-1}}}
  \sum_{m_1>\dotsb>m_{r-1}>m_r}
   \frac{(-1)^{d_1(m_1-m_2)+d_2(m_2-m_3)+\dotsb+d_{r-1}(m_{r-1}-m_r)+d_rm_r}}
    {m_1^{k_1}m_2^{k_2}\dotsb m_{r-1}^{k_{r-1}}m_r^{k_r}}\\
 =
  \sum_{m_1>\dotsb>m_{r'-1}>m_{r'}}
   \frac{(\frac12)^{-e_1(m_1-m_2)-e_2(m_2-m_3)-\dots-e_{r'-1}(m_{r'-1}-m_{r'})-e_{r'}m_{r'}}}
    {m_1^{k'_1}m_2^{k'_2}\dotsb m_{r'-1}^{k'_{r'-1}}m_{r'}^{k'_{r'}}},
  \label{eq:alt dual}
 \end{multline}
 where $|w_{01}|_{x_{-1}}$ is
 the total number of factors of $x_{-1}$ in $w_{01}$.
 For example, the duality formula w.r.t.\ $x_{-1}$ is
 \begin{align}
  -(-\log(2))
  =
  -\sum_{m=1}^\infty \frac{(-)^m}{m}
  =
  \sum_{m=1}^\infty \frac{1}{m}\left(\frac{1}{2}\right)^m
  =
  -\log\left(1-\frac{1}{2}\right),
  \label{eq:ex1}
 \end{align}
 and the duality formula w.r.t.\ $x_0x_1x_{-1}$ is
 \begin{align}
  -\sum_{m_1>m_2} \frac{(-)^{m_2}}{m_1^2 m_2}
  =
  \sum_{m_1>m_2} \frac{1}{m_1m_2^2}\left(\frac{1}{2}\right)^{m_1-m_2}.
    \label{eq:ex2}
 \end{align}
 Equations \eqref{eq:ex1} and \eqref{eq:ex2} are respectively
 the duality formulas of Examples 6.8 and 6.1 of \cite{MR1709772}.
 The left hand side of \eqref{eq:alt dual} is the multiple $L$-value of
 modulus $2$ but unfortunately the right hand side is not.
\end{exam}

\section{Duality from the symmetry of $\sing$}\label{sec: dual}

For singular points $\sing$,
we consider the subgroup $G_\sing$ of the linear transforms PSL$(2;\C)$
preserving $\sing\cup\{\infty\}$.
Defining the action of $\sigma\in G_\sing$ on the $x_c\in\A_\sing$ as
\begin{align*}
 \sigma(x_c) = x_{\sigma(c)} - x_{\sigma(\infty)}
\end{align*}
and extending as the homomorphism w.r.t.\ ordinary product
(so $\sigma$ is a $\sh$-homomorphism),
$G_\sing$ acts on $\A_\sing$.
(Here $x_\infty$ means $0$.)
Then $G_\sing$ acts on $\Liz{\sing}{ab}{w}{z}$ on the left and right by
\begin{align*}
 \left(\sigma\Li{\sing}{ab}\right)\z{w}{z}
  &= \Liz{\sing}{ab}{\sigma^{-1}(w)}{z},\\
 \left(\Li{\sing}{ab}\sigma\right)\left(w;z\right)
  &= \Liz{\sing}{ab}{w}{\sigma(z)}.
\end{align*}

\begin{prop}
For $\sigma, \tau \in G_\sing$, and $c\in\sing$ we have
 \begin{multline*}
 \frac{d}{dz}\left( \tau\Li{\sing}{ab}\sigma \right)\z{x_cw}{z}\\
 =
 \left\{
  \frac1{z-(\tau\circ\sigma)^{-1}(c)}-
  \frac1{z-(\tau\circ\sigma)^{-1}(\infty)}
 \right\}
 \left( \tau\Li{\sing}{ab}\sigma \right)\z{w}{z}.
 \end{multline*}
\end{prop}
\begin{proof}
If $\tau(\infty) \not\in \{c, \infty\}$
 \begin{align*}
 &\frac{d}{dz}\left( \tau\Li{\sing}{ab}\sigma \right)\z{x_cw}{z}
 =
 \frac{d}{dz} \Liz{\sing}{ab}
  {(x_{\tau^{-1}(c)}-x_{\tau^{-1}(\infty)})\tau^{-1}(w)}{\sigma(z)}\\
  &\quad=
  \begin{aligned}[t]
  &\left\{
   \frac1{z-\sigma^{-1}(\tau^{-1}(c))} - \frac1{z-\sigma^{-1}(\infty)}
  \right\}\left( \tau\Li{\sing}{ab}\sigma \right)\z{x_cw}{z}\\
  &\qquad-
  \left\{
   \frac1{z-\sigma^{-1}(\tau^{-1}(\infty))} - \frac1{z-\sigma^{-1}(\infty)}
  \right\}
  \left( \tau\Li{\sing}{ab}\sigma \right)\z{x_cw}{z}
  \end{aligned}\\
 &\quad=
 \left\{
  \frac1{z-(\tau\circ\sigma)^{-1}(c)}-
  \frac1{z-(\tau\circ\sigma)^{-1}(\infty)}
 \right\}
 \left( \tau\Li{\sing}{ab}\sigma \right)\z{w}{z}.
 \end{align*}
 If $\tau(\infty)= c$,
 \begin{align*}
 &\frac{d}{dz}\left( \tau\Li{\sing}{ab}\sigma \right)\z{x_cw}{z}
 =
 -\frac{d}{dz} \Liz{\sing}{ab}
  {x_{\tau^{-1}(\infty)}\tau^{-1}(w)}{\sigma(z)}\\
  &\quad=
  \begin{aligned}[t]
  -
  \left\{
   \frac1{z-\sigma^{-1}(\tau^{-1}(\infty))} - \frac1{z-\sigma^{-1}(\infty)}
  \right\}
  \left( \tau\Li{\sing}{ab}\sigma \right)\z{x_cw}{z}
  \end{aligned}\\
 &\quad=
 \left\{
  \frac1{z-(\tau\circ\sigma)^{-1}(c)}-
  \frac1{z-(\tau\circ\sigma)^{-1}(\infty)}
 \right\}
 \left( \tau\Li{\sing}{ab}\sigma \right)\z{w}{z}.
 \end{align*}
 In the same way the case $\tau(\infty)=\infty$ can be proved.
\end{proof}

In particular take $\tau = \sigma^{-1}$ and we have
\begin{align*}
 \frac{d}{dz}\left( \sigma^{-1}\Li{\sing}{ab}\sigma \right)\z{x_cw}{z}
 =
 \frac{1}{z-c} \left( \sigma^{-1}\Li{\sing}{ab}\sigma \right)\z{w}{z},
\end{align*}
i.e.\ 
\begin{align*}
 (\sigma^{-1}H_\sing^{ab}\sigma)(z)
 =
 \sum_W \left(\sigma^{-1}\Li{\sing}{ab}\sigma\right)\z{w}{z} W
\end{align*}
becomes the solution of \eqref{eq:ext KZ} again.
To determine this function,
we need to know the asymptotic property around some point.

$G_\sing$ also acts on $\ncfps{X_c;c\in\sing}$ as
\begin{align*}
 \sigma(X_c) = X_{\sigma(c)},
\end{align*}
where $X_\infty=-\sum_{c\in\sing}X_c$.
Then
\begin{align*}
 \sum_W \sigma(w)\ W
 &=
 \sum_W w\ \sigma^{-1}(W)\\
 &=
 \exp(x_a\sigma^{-1}(X_a))
 \left(\sum_W \reg^{ab}(w)\sigma^{-1}(W) \right)
 \exp(x_b\sigma^{-1}(X_b))
 .
\end{align*}
Applying $\Liz{\sing}{ab}{\cdot}{z}$ to this equation, we get
\begin{align*}
 &
 (\sigma^{-1}H_\sing^{ab}\sigma)(z)
 =
 \sum_W\Liz{\sing}{ab}{w}{\sigma(z)} \sigma^{-1}(W)\\
 &=
 \left(\frac{\sigma(z)-a}{b-a}\right)^{X_{\sigma^{-1}(a)}}
 \left(
  \sum_W \Liz{\sing}{ab}{\reg^{ab}(w)}{\sigma(z)} \sigma^{-1}(W)
 \right)
 \left(\frac{\sigma(z)-b}{a-b}\right)^{X_{\sigma^{-1}(b)}}.
\end{align*}
In particular assume that $\sigma(a)=b$ and $\sigma(b)=a$.
Then $\sigma$ is expressed as
\begin{align*}
 \sigma(z) =
 \begin{cases}
  a + b - z& \text{if $\sigma(z)=\infty$},\\
  \displaystyle
   \alpha - \frac{(\alpha - a)(\alpha - b)}{\alpha - z}
  \quad\text{for some $\alpha\in\C$}&
    \text{if $\sigma(z)\neq\infty$},
 \end{cases}
\end{align*}
and it is easy to check that $\sigma$ is involutive.
Here $\alpha$ is the image of $\infty$.
Under this assumption, we have
\begin{multline*}
 (\sigma^{-1}H_\sing^{ab}\sigma)(z)\\
 =
 \left(
  \frac{z-b}{a-b} \frac{a-\sigma(\infty)}{z-\sigma(\infty)}
 \right)^{X_b}
 \left(
  \sum_W \Liz{\sing}{ab}{\reg^{ab}(w)}{\sigma(z)} \sigma^{-1}(W)
 \right)\\
 \times
 \left(
  \frac{z-a}{b-a} \frac{b-\sigma(\infty)}{z-\sigma(\infty)}
 \right)^{X_a},
\end{multline*}
where
\begin{align*}
 \frac{a-\sigma(\infty)}{z-\sigma(\infty)}
 =
 \frac{b-\sigma(\infty)}{z-\sigma(\infty)}
 =1
 \qquad\text{if $\sigma(\infty)=\infty$}.
\end{align*}
Thus we have
\begin{align*}
 (\sigma^{-1}H_\sing^{ab}\sigma)(z)
 \left( \frac{z-a}{b-a} \right)^{-X_a}
 \longrightarrow
 \left( \frac{b-\sigma(\infty)}{a-\sigma(\infty)} \right)^{X_a}
 \qquad (z\to a),
\end{align*}
in other words
\begin{align*}
 (\sigma^{-1}H_\sing^{ab}\sigma)(z)
 \left( \frac{a-\sigma(\infty)}{b-\sigma(\infty)} \right)^{X_a}
 =
 H_\sing^{ba}(z).
\end{align*}
Moreover assume that
$a$ and $b$ are the nearest points to each other in $\sing$.
Then by the above expression of $H_\sing^{ba}(z)$,
we can compute the connection matrix as
\begin{align*}
 &\Phi_\sing^{ab}(X_c;c\in\sing)
 = H_\sing^{ba}(z)^{-1} H_\sing^{ab}(z)\\
 &=
 \left( \frac{b-\sigma(\infty)}{a-\sigma(\infty)} \right)^{X_a}
 \sum_W
  \left(
   \sum_{w_1w_2=w}
    \left(
     \sigma^{-1}\Li{\sing}{ab}\sigma \right) \z{S(w_1)}{z}
      \Liz{\sing}{ab}{w_2}{z}
    \right)
 W
 \\
\begin{split}
  &=
 \begin{aligned}[t]
 &\left(
  \frac{b-a}{z-a} \frac{z-\sigma(\infty)}{a-\sigma(\infty)}
 \right)^{X_a}
 \left(
  \sum_W \Liz{\sing}{ab}{\reg^{ab}\circ\sigma\circ S(w)}{\sigma(z)} W
 \right)\\
 &\hspace{7.5cm}\times
 \left(
  \frac{a-b}{z-b} \frac{z-\sigma(\infty)}{a-\sigma(\infty)}
 \right)^{X_b}
 \end{aligned} 
 \\ 
 &\qquad\times
  \left(\frac{z-a}{b-a}\right)^{X_a}
   \left(\sum_W\Liz{\sing}{ab}{\reg^{ab}(w)}{z}W\right)
    \left(\frac{z-b}{a-b}\right)^{X_b}.
\end{split}
\end{align*}
In the above equations
by tending $z$ to $a$ and $b$,
we have the next theorem.

\begin{thm}
 For a linear transform $\sigma$
 which preserves $\sing\cup\{\infty\}$
 and interchanges $a$ and $b$,
 two points nearest each other in $\sing$,
 we have
 \begin{multline*}
 \Phi_\sing^{ab}(X_c;c\in\sing)
 \\
 =
 \left(\frac{b-\sigma(\infty)}{a-\sigma(\infty)}\right)^{X_a}
 \left(
  \sum_W \Lv_\sing^{ab}(\reg^{ab}\circ\tau(w)) W
 \right)
 \left(\frac{b-\sigma(\infty)}{a-\sigma(\infty)}\right)^{X_b}
 ,
 \end{multline*}
 where $\tau:=\sigma\circ S$.
 Equivalently,
 the inverse of $\Phi_\sing^{ab}(X_c;c\in\sing)$ can be written as
 \begin{align*}
  \Phi_\sing^{ab}(X_c;c\in\sing)^{-1}
  =
  \left(\frac{a-\sigma(\infty)}{b-\sigma(\infty)}\right)^{X_b}
  \Phi_\sing^{ab}(\sigma(X_c);c\in\sing)
  \left(\frac{a-\sigma(\infty)}{b-\sigma(\infty)}\right)^{X_a}.
 \end{align*}
 In the special case, for $w\in\A_\sing^{ab}$
 there holds Euler's inversion formula
 \begin{align*}
 \sum_{w_1w_2=w}
  \Liz{\sing}{ab}{\tau(w_1)}{\sigma(z)}\Liz{\sing}{ab}{w_2}{z}
 = \Lv(w)
 \end{align*}
 and a ``duality formula''
 \begin{align}
 \Lv(w)=\Lv(\tau(w)).
 \label{eq: involutive dual}
 \end{align}
\end{thm}

The map $\tau$ is involutive
because $\sigma\circ S = S\circ\sigma$
and both $\sigma$ and $S$ are involutive.
So it is appropriate to call the equation
\eqref{eq: involutive dual} the ``duality formula''.

\begin{exam}
 Let $\sing=\{0,1\}$, $\sigma(z)=1-z$ and $b=0$ (i.e.\ $a=1$).
 Then $\tau$ is the anti-involution defined by
 \begin{align*}
  x_0 &\mapsto -x_1,& x_1\mapsto& -x_0.
 \end{align*}
 Define $\h^0 = \A_\sing^{10}$, $x=x_0$, $y=-x_1$,
 and $\zeta:\h^0\to\R$ for $w=x^{k_1-1}y \dotsb x^{k_r-1}y\in\h^0$ by
 \begin{align*}
  \zeta(w) := \Lv_\sing^{ab}(w)
  &= \sum_{m_1>m_2>\dots>m_r} \frac{1}{m_1^{k_1}m_2^{k_2}\dots m_r^{k_r}}\\
  &= \zeta(k_1,k_2,\dots,k_r);
 \end{align*}
 then $\tau$ is the anti-involutive mapping
 \begin{align*}
  x &\mapsto y,& y\mapsto& x,
 \end{align*}
 satisfying the next relation;
 \begin{align*}
  \zeta(\tau(w)) = \zeta(w).
 \end{align*}
 This is the formulation of the duality formula by \cite{MR99e:11119}.

 The inverse of $\Phi_\sing^{10}(X_0,X_1)$ can be computed as
 \begin{align*}
  \Phi_\sing^{10}(X_0,X_1)^{-1}
  &= \left(\sum_W \Lv_\sing^{10}(\reg^{10}\circ\sigma\circ S(w))W\right)^{-1}\\
  &= \sum_W \Lv_\sing^{10}(\reg^{10}\circ\sigma(w))W
   = \sum_W \Lv_\sing^{10}(\reg^{10}(w))\sigma(W)\\
  &= \Phi_\sing^{10}(\sigma(X_0),\sigma(X_1)) = \Phi_\sing^{10}(X_1,X_0).
 \end{align*}
 This relation is nothing but \eqref{eq: associator dual}.
\end{exam}

\begin{exam}
 Take $\sing=\{0,\pm1,\pm\imu\}$, $\sigma=(1-z)/(1+z)$ and $b=0$.
 Then for any word of $\A_\sing^{ab}$,
 $\Lv_\sing^{ab}$ corresponds to the multiple $L$-value of $m=4$.
 For $w=x_0^{k_1-1}x_{c_1}\dotsm x_0^{k_{r-1}-1}x_{c_{r-1}}x_0^{k_r-1}x_{c_r}$
 the corresponding value is
 \begin{align*}
  \Lv_{\sing}^{ab}(w)
  &=
  (-1)^r
  \sum_{m_1=r}^\infty
   \sum_{m_1>\dots>m_{r-1}>m_r>0}
    \frac{c_1^{-m_1+m_2} \dotsm c_{r-1}^{-m_{r-1}+m_r} c_r^{-m_r}}
     {m_1^{k_1} \dotsm m_{r-1}^{k_{r-1}} m_r^{k_r}}\\
  &=
   (-1)^r
   L(k_1,\dots,k_{r-1},k_r;a_1,\dots,a_{r-1},a_r),
 \end{align*}
 where $(\imu)^{-a_i}=c_i$.
 The involution $\sigma$ interchanges the singular points
 \begin{align*}
  \sigma&:&
     0 &\longleftrightarrow 1,&
  \imu &\longleftrightarrow -\imu,&
    -1 &\longleftrightarrow \infty,
 \end{align*}
 so the action of $\sigma$ on $\A_\sing$,
 the action of $\sigma$ on $\ncfps{X_c;c\in\sing}$,
 and
 the action of the anti-involution $\tau=\sigma\circ S$ on $\A_\sing$
 are as follows:
 \begin{align*}
  \begin{aligned}
  \sigma&:&
  x_0 &\mapsto x_1 - x_{-1},&
  x_1 &\mapsto x_0 - x_{-1},&
  x_{-1} &\mapsto - x_{-1},\\
  && 
  x_{\scriptscriptstyle \imu} &\mapsto x_{\scriptscriptstyle -\imu} - x_{-1},&
  x_{\scriptscriptstyle -\imu} &\mapsto x_{\scriptscriptstyle \imu} - x_{-1},&&
  \end{aligned} 
 \end{align*}
 \begin{align*}
 \sigma&:&
  X_0 &\mapsto X_1,&
  X_1 &\mapsto X_0,&
  X_{-1} &\mapsto -X_0-X_1-X_{\imu}-X_{-1}-X_{-\imu},
  \notag\\
  &&
  X_{\scriptscriptstyle \imu} &\mapsto X_{\scriptscriptstyle -\imu},&
  X_{\scriptscriptstyle -\imu} &\mapsto X_{\scriptscriptstyle \imu},&&
 \end{align*}
 \begin{align*}
  \begin{aligned}
   \tau&:&
  x_0 &\mapsto - x_1 + x_{-1},&
  x_1 &\mapsto - x_0 + x_{-1},&
  x_{-1} &\mapsto x_{-1},\\
  && 
  x_{\scriptscriptstyle \imu} &\mapsto - x_{\scriptscriptstyle -\imu} + x_{-1},&
  x_{\scriptscriptstyle -\imu} &\mapsto - x_{\scriptscriptstyle \imu} + x_{-1}.&
  &
  \end{aligned}
 \end{align*}
 We have the duality formula of multiple $L$-values
 \begin{align*}
  \Lv_\sing^{ab}(\tau(w)) = \Lv_\sing^{ab}(w),\qquad w\in\A_\sing^{ab},
 \end{align*}
 which is equivalent to the formula:
 \begin{align*}
  \begin{split}
   &\Phi_\sing^{ab}(X_0,X_1,X_{\imu},X_{-1},X_{-\imu})^{-1}\\
  &=
   2^{X_0}
   \Phi_\sing^{ab}(\sigma(X_0),\sigma(X_1),
     \sigma(X_{\imu}),\sigma(X_{-1}),\sigma(X_{-\imu}))
   2^{X_1}.
  \end{split} 
 \end{align*}
 The transform $\sigma$ also preserves the subset
 $\{0,1,-1,\infty\}\subset\sing\cup\{\infty\}$.
 So the above duality formula can be restricted
 from $\sing$ to $\sing'=\{0,1,-1\}$.
 The duality formula for $\sing'$ is
 the reproduction of equation (127) of \cite{MR2002a:81180}.
\end{exam}

\section*{Acknowledgements} 
 The author would like to express his grantitude to
 Saburo Kakei, Michitomo Nishizawa and Yuji Yamada
 for useful communications.
 He also thanks Kimio Ueno for valuable advice and comments.

\def\cprime{$'$} \def\cprime{$'$}

\begin{tabular}{l}
 {Department~of~Mathematical~Sciences},\\
 School of Science and Engineering,\\
 WASEDA UNIVERSITY,\\
 3-4-1, Okubo Shinjuku-ku,\\
 Tokyo 169-8555, Japan\\
 {\ttfamily okuda@gm.math.waseda.ac.jp}
\end{tabular}

\end{document}